\documentclass[12pt]{article}
\usepackage{amssymb}

\newtheorem{theorem}{Theorem}

\newtheorem{claim}[theorem]{Claim}

\newtheorem{corollary}[theorem]{Corollary}

\newtheorem{definition}[theorem]{Definition}

\newtheorem{lemma}[theorem]{Lemma}

\newtheorem{remark}[theorem]{Remark}

\newcommand{\be}{\begin{equation}}
\newcommand{\ee}{\end{equation}}
\newcommand{\ba}{\begin{eqnarray}}
\newcommand{\ea}{\end{eqnarray}}
\newcommand{\ban}{\begin{eqnarray*}}
\newcommand{\ean}{\end{eqnarray*}}
\newcommand{\Tr}{{\rm{Tr }}}
\newcommand{\bcw}{\mathbin{\bigcirc\mkern-15mu\wedge}}
\begin{document}
\title{A class of variational functionals in conformal geometry}

\author{Sun-Yung Alice Chang \\ Princeton University\thanks{The first author is partially supported by NSF grant DMS-0245266. }
 \and Hao Fang \\University of Iowa \thanks
{The second author is partially supported by NSF grant DMS-0606721.}
}


\maketitle

\abstract{We derive a class of variational functionals which arise naturally in
conformal geometry. In the special case when the Riemannian manifold is locally
conformal flat, the functional coincides with the well studied functional which is
the integration over the manifold of the k-symmetric function of the Schouten tensor of the metric on the manifold.}

\section{Introduction}
The purpose of this article is to derive a class of variational
functionals which arise naturally in conformal geometry. Recall on
a Riemannian manifold $(M^n, g)$, $n \geq 3$, the full Riemannian curvature tensor $Rm$ decomposes as
 $$
Rm=W \,\, \oplus \,\,P \bcw  g,
$$
where $W$ denotes the Weyl tensor,
$$
P \,= \, \frac {1}{n-2} ( Ric - \frac{R}{2(n-1)}g)
$$
denotes the Schouten tensor, and  $\bcw$ is the Kulkarni-Nomizu
wedge product (see [Be87, pp. 110]). Under a conformal change
of metrics $g_w = e^{2w} g $, where $w$ is a smooth function
over the manifold, the Weyl curvature changes
pointwisely as  $W_{g_w} = e^{2w} W_g$. Thus, all the information of the
Riemannian tensor under a conformal change of metrics is reflected
by the change of the Schouten tensor, and it is natural to study
the elementary symmetric function $ \sigma_k (g^{-1}P_g)$ (which
we later denote as $\sigma_k (g)$) of the eigenvalues of the
Schouten under the conformal change of metrics. For example when
$k =1$, $\sigma_1 (g) = \frac {1}{2(n-1)} R_g$, where $R_g$
denotes the scalar curvature of g. The study of the equation
$\sigma_1 (g) = $ constant under conformal change of metrics is
the classical Yamabe problem. In ~{[V00]}, Viaclovsky proved the
following statements: Consider the functional
\begin{equation}
 F_k (g) : = \int_{M^n} \sigma_k (g) dv_g,
\end{equation}
\vskip .1in (i): When $k=1 \,\, { or} \,\,\, 2$, and  $ 2k <n$,
$F_k$ is variational in the conformal class of metrics $g_w \in
[g]$ with fixed volume one, i.e. the extremal metric
for the functional in this class of metrics, when achieved,
satisfies the equation
$$
 \sigma_{k}{(g_w)} = {\rm constant}.
$$
\vskip .1in (ii) When $k \geq 3$ and $ 2k <n$, assertion in (i)
only holds when the manifold $M^n$ is locally conformally flat.
\vskip .1in (iii) When $ k=2, \, n=4$, $F_2 (g)$ is conformally
invariant; while for $ k = \frac n 2$ and $ k \geq3 $, $F_k (g) $
is conformally invariant only when the manifold $(M^n,g)$ is
locally conformally flat. \vskip .1in

We remark that in ~{[BG06]}, Branson and Gover have also proved that
the metric being locally comformally flat is also a necessary 
condition for the statement of (ii) above.

In this article, we generalize the role played by the
curvature polynomial $\sigma_k (g) $  to a new class of curvature
invariant, $v^{(2k)}(g)$, so that: $v^ {(2k)} (g) $ agrees with
$\sigma_k(g)$ when $(M^n, g)$ is locally conformally flat;  the
functional
\begin{equation}
{\cal F}_{k}(g) = \int_{M^n} v^{(2k)}(g) dv_g
\end{equation}
 satisfies the variational property in the statement (ii) above for all $ 2k <n$;
$ {\cal F}_{k} (g)$ is conformally invariant when $ 2k =n$ for all
Riemannian manifolds $(M^n, g)$. \vskip .1in

Our construction is closely related to the recent construction for
the $Q$ curvature of R. Graham and Juhl~[GJ06]. To state our
result, we first recall some definitions and basic results.
\vskip.1in
\begin{definition} Given $(X^{n+1},M^n, g^+)$ with smooth boundary $\partial X = M^n$.
Let $r$ be a defining function for $M^n$ in $X^{n+1}$ so that
$r>0$ in $X$ while $r=0 \,\, {on} \; M \;{ and} \,\,  dr|_{M}\neq 0
$. We say $g^+$ is conformally compact, if there exists such $r$
so that $(X^{n+1}, r^2 g^+)$ is a compact manifold. We say
$(X^{n+1}, M^n, g^+)$ is  conformally compact Einstein if $g^+$ is
Einstein (i.e. $ Ric_{g^+} = c g^+ $ for some constant c ). We
call $ g^+ $ Poincare metric if $ Ric_{g^+} = - {n} g^+ $.
\end{definition}

Fefferman-Graham~[FG85] proved that for any given $(M^n,
g_0)$, there is an extension,  $g^{+}$, which is "asymptotically
Poincare Einstein" in a neighborhood of $ M^n$, i.e. on $ [0,
\epsilon ) \times M^n $ for some positive $\epsilon$.

For $c$ sufficiently small, we have $X_{c} \,= \, \{ r\leq
c\}\subset X$ is diffeomorphic to $[0,c] \times M $. Hence, given
a local coordinate chart $(x^{1},\cdots, x^{n})$ on $M$,
$(r,x^{1},\cdots, x^{n})$ is a coordinate chart on $X_{c}$. Thus,
we can write

 \begin{equation}
 g_{+}={1\over r^{2}}(dr^{2}+h_{ij}(r,x)dx^{i}dx^{j}),
 \end{equation}
 where $h(r,\cdot)$ is a metric defined for $M_{c}:={r=c}\subset X$. We
remark that
in this expression, $(X, r^2 g_{+})$ is a compact manifold.  

 \[
 dvol_{h}(r,x)=\sqrt{\det h_{ij}(r,x)} dx^{1}\cdots dx^{n}.
 \]

 We now expand the quantity $\sqrt{{\det h_{ij}(r,x)}\over{ \det h_{ij}(0,x)}}$  in an expansion near $r=0$
as
 \begin{equation}
 \sqrt{{\det h_{ij}(r,x)}\over{ \det h_{ij}(0,x)}}
 =\sum_{k=0}^{\infty}{v^{(k)}(x,h_{0})r^{k}},
 \end{equation}
 where $v^{(k)}(x,h_{0})$ is a curvature invariant of the metric 
$h_{0}:={ h_{ij}} {(0,\cdot )} $.

In~{[G00]}, R. Graham asserted that $v^{(k)}$ vanished for $k$
odd and $ 2k < n$; furthermore, he has established that
when $n$ is even, the quantity $\int_{M}v^{(n)}dvol_{h_{0}}$
is conformally invariant over the conformal class of metrics of
$[h_{0}]$. This quantity is related to the conformal invariant term in
the expansion of ``renormalized volume'' in the study of conformal
field theory. In the later papers, {[GZ03]}, {[FG02]} and
{[GJ06]}, the authors have also established that the quantity
is the same of the integral of the $Q$ curvature and is an
important global conformally invariant term.

Another fact which has been pointed out in [GJ06], based on a result
in {[SS00]},  is that
when $h_0$ is locally conformally flat, then
 \begin{equation}
 v^{(2l)} (h_0) \, = \, (-2)^{-l}\sigma_{l}(h_0).
 \end{equation}

 In this article, we prove the following

 \begin{theorem}\label{main}
 For any metric $h$ on $M$ and $ 2k\leq n=\dim M$, define the functional
 \begin{equation}
 {\cal F}_{k}(h)=\int_{M}{v^{(2k)}(h)dvol_{h}/ (\int_{N}dvol_{h})^{{n-2k}\over n}},
 \end{equation}
 then ${\cal F}_{k}$ is variational within the conformal class when $ 2k
< n $; i.e., the critical metric in $[h]$ satisfies the equation
 $$ v^{(2k)}={\rm constant}.
$$
 For $n=2k$,  $F_{\frac n2}(h)$ is constant in the conformal class $[h]$.
 \end{theorem}

 \begin{remark}
 As we have mentioned before, the case $n=2k$ in the theorem above
 has been established earlier in~{[G00]}.
 \end{remark}

 For $k=1,\, 2\, $ cases, the new curvature invariant $v^{(2k)}$  turns out to agree, up to a scale, with
the well-studied curvature polynomial $\sigma_k (g)$. Actually we
have
 $$v^{(2)} (g) \, =\, - \frac 12 \sigma_{1}(g),$$
 $$ v^{(4)} \, = \, \frac 14 \sigma_{2}(g).$$
 For $k=3$, Graham and Juhl ({[GJ06]}, page 5) have also listed the following formula for $v^{(6)}$
 
 $${ v}^{(6)}(g) \, = \, - \frac {1}{8} \left[\sigma_{3}(g) \, + \,  {1\over{3(n-4)}}(P_g)^{ij}{(B_g)}_{ij}\right],$$
 where  ${(B_g)}_{ij} : ={1\over{n-3}} \nabla^k \nabla^l W_{likj} +{1\over{n-2}} 
R^{kl}W_{likj}$ is the Bach tensor of the metric.
 
In this article, we carry out the 
computation for $v^{(2)}$ and $v^{(4)}$. As the computation indicates, a straight
forward computation of $v^{(6)}$  is
quite complicated. Instead,   we derive the variational
properties of $v^{(6)}$ under conformal change of metrics directly, which is another verification of our main theorem in this special case.
Another purpose of the derivation is to derive the variational formulas of $v^{(6)}$ under
conformal change of metrics. We believe the study of the PDE $
v^{(6)}(g_w) = constant $ will
be of interest to problems in conformal geometry.
Another interesting question is the ``uniqueness''
problem of curvature invariants which are extensions of $\sigma_{k}(g)$ invariants in the locally conformally flat case and satisfy the
properties as $v^{(2k)}$ in Theorem 2 above. We hope to address these
two problems in a future work.

This article is organized as follows: In Section 2, we prove Theorem~\ref{main}; In Section 3,
we explicitly compute $v^{(2)}$ and $v^{(4)}$ for all dimension $n$; In Section 4, we show the variational property of ${v}^{(6)}$ and discuss some properties 
of this curvature invariants under conformal change of metrics.

\vskip .2in

\section{Proof of Theorem~\ref{main}}

Suppose $(X^{n+1}, M^n, g_{+})$ is a conformally compact Einstein
manifold. Let $r$ be an arbitrary smooth defining function
for $M= \partial X$ defined near $M$ and set $\bar g= r^2 g+$.
We now recall some basic properties of the metric $g_{+}$
in this setting with respect to the changing of defining functions.

\begin{lemma} ([G00], Lemma 2.1 and 2.2]):\label{lemma2.1}\newline

(a) A metric on M in the conformal infinity of $g^+$ determines a
unique defining function $r$ in a neighborhood of $M$ such that
$\bar g |_{TM} $ is the prescribed boundary metric with
$|dr|_{\bar g}^2 =1$.

(b) \be
 g_{+}  = r^{-2} (dr^2 + g_r)
 \label{2.1} \ee on
$[0,\epsilon) \times M$ for some $\epsilon >0$. Furthermore, \be
g_r=g + r^2 g^{(2)} + r^4 g^{(4)} + ..+r^{n} g^{(n)} + \tilde h
r^n log r + \cdots \label{2.2} \ee when n is even, with $g = r^2
g^+ |_{M}$, and symmetric tensors $g^{(2)}$, $g^{(4)}$, ... up to $g^{2
(n-1)}$ and $Tr_gg^{(n)}$ are determined by $g$, 
and $Tr_g \tilde h = 0$.

(c) Let $r$ and $\hat r$ be two special defining functions as in (a) associated
with two different conformal representatives in the conformal class of metrics
in $[g]$; then
\be \hat r = r e^w \label{2.3} \ee
for a function $w$ on $ [0, \epsilon) \times M$ satisfying
\be
w_{r}+r(w_{r}^{2}+|d_M w |^2)= 0.
\label{2.4}
\ee
Furthermore, the power series expansion of $w$ at $r=0$ consists only of even power of
$r$ up through and including the $r^{n+1}$ term.
 \end{lemma}

Now for a fixed smooth function $\phi$ defined on $M$, we consider a family of
conformal metrics $g_t = e^{2t \phi} g$  on $M$.  By Lemma~\ref{lemma2.1}   there exist functions $r_{t}$ on a neighborhood of $M$ in $X$ such that
\begin{equation}\label{g+}
 g_{+}={1\over r_{t}^{2}}(dr_{t}^{2}+h_{t}(r_{t},\cdot)),
 \end{equation}
 where $h_{t}(c,\cdot)$ is a metric defined on $M_{t,c}=\{r_{t}=c\}\subset X$. Furthermore,
 we have the following asymptotic expansion
\be h_t =  g_t + r_t^2 g_t^{(2)} + \cdots. \label{2.5}
\ee

For a point $p\in X$, define  $$w(t,p)=\log ({r_{t}(p)\over r(p)}).$$
On the boundary, we have that $w(t,\cdot)|_{M}=t\phi(\cdot)$. Thus,
$w(t,p)=w(t,r,x)$ is a smooth extension of $t \phi(x) $ to $X_{t,c}:=\{r_{t}<c\}\subset X$ for some proper $c$.
 We then get the following from Lemma~\ref{lemma2.1}:

 \begin{corollary}  For $(t,x, r) \in [0,1] \times
[0, \epsilon) \times M$, we have
\be
\frac{\partial}{\partial r }|_{r=0} w(t,r,x) = 0.
\label{2.6}
\ee
\be
\frac{ \partial^2}{ \partial r^2} |_{r=0} w(t,r,x)
= - \frac {1}{2} t^2 |d_M \phi (x)|_{g} ^2 .
\label{2.7}
\ee
\be
\frac{\partial^k}{\partial r^k}|_{r=0} w(t,r,x)
 = 0, \rm {for \,\, k \,\, odd}.
\label{2.8}
\ee
\be
\frac{\partial^k}{\partial r^k}|_{r=0} w(t,r,x)
 = O (t^2),
  \rm {for \, \,  k \,\, even, \,\, 0 <  k  < n}.
\label{2.9}
\ee
In particular, we have
\be
w(t,r,x)= t\phi(x) -\frac{1}{4} t^{2}|\nabla_{M}\phi(x)|^{2}r^{2}+ O(t^2 r^4).
\label{2.10}
\ee
and
\be
\frac {d}{dt}|_{t=0} w(t,r,x) = \phi (x),
\label{2.11}
\ee
independent of the choice of the defining function $r$.
\end{corollary}

For future use, we define a useful vector field associated to the conformal metric variation.

First we notice that, fix an $\epsilon$ small enough, for a given point $p \in [0, \epsilon) \times M= X_{0,\epsilon}$,
for each t, we can assign a local coordinate chart $ p = (r_t, x_t)\in X_{t, \epsilon}$ for $t\in[0,1]$,
with $ r_0 =r, x_0 = x$ and with $ x_t$ defined as $pr_{t}(p)$, the projection image of $p$ onto $M$ under the metric $r_{t}^{2}g_{+}$.

For each fixed $c < \epsilon$, denote $M_{c,t} = \{ p \in X| r_t
=c \} $, then by Lemma~\ref{lemma2.1}, the set $M_{c,t}$ is
diffeomorphic  to  $M$ via the projection with respect to the
metric $r_{t}^{2}g_{+}$. Hence, they are diffeomorphic to each
other. For a fixed level set $M_c := M_{c,0}$, since the projectoin $pr_{t}$ is a small perturbation of $pr_{0}$ when $t$ is small, the map $ p \to (c,
x_t)$ gives arise to a diffeomorphism of $M_c$; hence it
introduces a vector field \be\label{vf} F_c := F_{c,0} = \frac
{d}{dt}|_{t=0} (c, x_t) \ee on $TM_c$. Since $M_{c}$ is  naturally
diffeomorphic to $M$, without confusion, we also denote the
pushed-forward vector field on $M$ as $F_{c}$.

It is worth pointing out that: the vector field is induced from the one parameter family of points $x_{t}$, which depends on the original point $p$; hence the induced vector field  depends on the choice of  $c$.

Since we have a family of diffeomorphisms to identify a
neighborhood of $M$ in $X$ with $[0, \epsilon) \times M $, given a
local coordinate chart $(x^{1},\cdots, x^{n})$ on $M$,
$(r_{t},x^{1},\cdots, x^{n})$ is a coordinate chart on $M_{c,t}$,
for each $t$. Thus, a given point $p\in X$ near $M$ can be
represented in these coordinate system as $(r_{t},x_{t})$,
respectively.

We now consider the volume of $g_{+}$ at a given point $p\in M$. For future convenience, we will omit the foot index $0$, and denote $h=h_0$, $g= g_0 = h_0|_{TM} $,  $r=r_0$, $x=x_0$ etc. Thus, $p=(r,x)$ and the metric $r^{2}g_{+}$ is compact.

Recall $r_t (p) = r e^{ w (t, p)} $ near $M$. For each $t$, by~(\ref{g+}) and~(\ref{2.5}),
\begin{equation}\label{volume}
dvol_{g^{+}}(p)=r_{t}^{-n-1}dr_{t} dvol_{h_{t}}(r_t, x_{t})
={r_{t}^{-n-1} \sqrt{{\det h_{t}(r_{t},x_{t})}\over{ \det g_t (x_{t})}}dr_{t}dvol_{g_{t}}(x_{t})},
\end{equation}
and
$$
\sqrt{{\det h_{t}(r_{t},x_{t})}\over{ \det
g_{t}(x_{t})}}=\sum_{k}v^{(k)}(x_{t},g_{t})r_{t}^{k}.
$$

Notice that, via diffeomorphism, we can view $dvol_{g_{t}}(x_{t})$
as a $n$-form on $M_{r_{t},t}$.

We now proceed by take the time derivative of Equation (\ref{volume}).
 For notational convenience, we define the following linear operator
 $$D={d\over dt}|_{t=0}.$$

 We prove a technical lemma. 

\begin{lemma} \label{lemma2.2} At the point $p=(r,x)$, the following formula hold:

(a)  $$D(dr_{t})(r,x)\wedge dvol_{g}(x)=\phi (x) dr\wedge dvol_{g}(x).$$

(b) $$D[v^{(k)}(x_{t},g_{t})]=(F_{r}v^{(k)})(x, g )+ D[v^{(k)}(x,
g_{t})],$$where the definition of $F_{r}$ is given in (\ref{vf})
and the remarks following (\ref{vf});

(c) $$D[dvol_{g_{t}}(x_{t})]=[{\cal L}_{F_{r}}(dvol_{g})](x)+n\phi(x)
dvol_{g}(x),$$
where $\cal L$ is the Lie derivative on $M_{r}$ with respect to the given vector field. \end{lemma}

\bigskip
\noindent{\bf Proof of Lemma~\ref{lemma2.2}.} Since
$r_{t}=e^{w}r_{0}=e^{w}r$, we have
$$dr_{t}= r_{t} dw \, + \, e^{w} dr.$$
Apply $D$ to both sides of above equation, and use
Lemma~\ref{lemma2.1} and ~(\ref{2.11}), \ban
D(dr_{t}) &=&  r d\phi \, + \, r \phi \, dw|_{t=0} \,  + \phi \, dr\nonumber \nonumber \\
&=& r d\phi \, + \, \phi dr. \ean We then wedge above expression
with $dvol_g (x)$, and observe that the term $ r d\phi (x) \wedge
dvol_{g}(x)=0$ at the point $p=(r,x)$,
 we have thus established statement (a).

Statement (b) follows directly by Leibniz rule and (\ref{vf}).

To prove statement (c), we apply Leibniz rule again and get
$$
D[dvol_{g_{t}}(x_{t})]=D[dvol_{g}(x_{t})]+D[dvol_{g_{t}}(x)].
$$
Notice that when restricted to $M$,
$g_{t}=e^{2t \phi}g $, the result follows easily.

\bigskip

We now give the proof of Theorem~\ref{main}.

\noindent{\bf Proof of Theorem~\ref{main}.}
 Starting with the following
basic equation

\ban
0 &=&  D(dvol_{g^{+}}(p))
 \nonumber \\
  &=&  D({r_{t}^{-n-1}(\sum_{k}v^{(k)}(x_{t},g_{t})r_{t}^{k})
dr_{t}dvol_{g_{t}}(x_{t})}),
\ean

we use Leibniz Rule and apply (17) and Lemma 5 above to get

\ba\label{new1} 0 &=& \sum_{k} r^{k-n-1} \{D[v^{(k)}(x_{t},g_{t})]
\, + \,
 k\phi (x) v^{(k)}(x,g) \}dr \wedge dvol_{g}(x)\nonumber \\
&+& \, \sum_{k} r^{k-n-1} dr {\cal L}_{F_{r}}  [v^{(k)}dvol_{g}(x)].
\ea

We integrate (~\ref{new1}) over $M_{r}$ which is identified to $M$ via
the canonical diffeomorphism. Since the form involving the Lie derivative is exact,  it will vanish after the integration. Thus, we get  \be dr \sum
r^{k-n-1}\int_{M} {\{D[v^{(k)}(x_{t},g_{t})]+k\phi (x)
v^{(k)}(x,g) \} dvol_{g}(x)}=0. \ee

Now notice that the above equation holds for all small $r$,
we prove the following identity:

\begin{claim}
\be \int_{M} {\{D[v^{(k)}(x, g_{t})]+k\phi (x) v^{(k)}(x,g) \}
dvol_{g}}=0. \ee
\end{claim}

We now finish the proof of Theorem~\ref{main}.
\bigskip

Given $g_{t}=e^{2t\phi}g $ a variation of metrics on $M$
 in the conformal class of $g$, denote $V=\int_{M}dvol_{g}$.
 Then, 
\ban D[{\cal F}_{k}(g_{t})] &=& {1\over V^{1- \frac
{2k}{n}}}\int_{M}D[v^{(2k)}(x, g_{t})]\
dvol_g (x) \, + \, n \int_{M}{v^{(2k)}(x,g) \phi (x) \ dvol_g (x)}  \nonumber \\
& &\ \ \ \ \ \ \ \ \ - {{n-2k}\over {n V^{ 2- \frac {n-2k}{n}}}} 
\int_{M}{n \phi  \ dvol_g}
\int_{M}{v^{(2k)}(x,g) \ dvol_g (x) } \nonumber \\
&=& {n-2k \over V^{1-\frac {2k}{n} }} \int_{M}[ v^{(2k)}-{{\int_{M} v^{(2k)}
(g,x) dvol_g (x) }\over V}] \phi  \, \ dvol_g. \ean It implies that,
when $n>2k$, the critical metric $g$ of  the functional ${\cal F}_{k}$
satisfies
$$
v^{(2k)}(g)=constant.
$$
When $n=2k$, the computation shows that the functional is invariant under the conformal deformation.

\bigskip\bigskip

\section{ Computation of $v^{(4)}$ for $n>4$}

In this Section, we verify the formula for $v^{(4)}$. In
particular, we prove that for any dimension $n > 4$, $v^{(4)}$
equals to $\sigma_{2}(A)$ up to a constant. We remark that this
formula is stated without proof in ~[GJ06], and the method of
derivation is known to experts in this field, thus we will be brief in our
derivation.

We start with the basic equations

\ba Ric (g_{+})&=&-n g_{+},\label{eqn1}  \\
g_{+}&=&{1\over r^{2}}(dr^{2}+h(r,\cdot))\label{eqn2}  \\
 h &=&  g + r^2 g^{(2)} +r^{4}g^{(4)}+ \cdots+ \tilde h r^n log
r+\cdots. \label{eqn3}\ea

For future convenience, we denote
\[
C_{k}=  g^{(k)}g^{-1}.
\]

First, we have the following Lemma, which follows from a simple
computation.
\begin{lemma}\label{v2v4}
\ba
v^{(2)}&=&{1\over 2}\Tr\ C_{2},\nonumber \\
v^{(4)}&=&{1\over 2}[\Tr C_{4}+\sigma_{2}(C_{2})-{1\over 4}(\Tr\
C_{2})^{2}].
\ea
\end{lemma}
We set up a local coordinate $\{x_{1},\cdots, x_{n}\}$ on $M_{r}$; thus,
$\{x_{1},\cdots, x_{n},x_{n+1}=r\}$ is a coordinate on $X$.
By (\ref{eqn1}) and (\ref{eqn2}), we get,

\be r
h''_{ij}+(1-n)h'_{ij}-h^{kl}h'_{kl}h_{ij}-rh^{kl}h'_{ik}h'_{lj}+{r\over
2}h^{kl}h'_{kl}h'_{ij}-2r\  Ric(h)_{ij}=0,\label{eqn4}\ee
where we use $'$ to denote the derivative with respect to $r$ and $Ric
(h)$ is the Ricci curvature of the submanifold $M_{r}$ with respect to
the restricted metric $h(r,\cdot)$.

We analyze (\ref{eqn4}) by (\ref{eqn3}). Since  $n > 4$ we get \ba
h'=2rg^{(2)}+4r^{3}g^{(4)}+O(r^{4}),\label{eq3}\nonumber \\
h''=2g^{(2)}+12r^{2}g^{(4)}+O(r^{3}). \label{eq4}
\ea
Studying the coefficient of $r$ in (\ref{eqn4}), we get a tensor
equation over $M$,
\be\label{eqn5}
2g^{(2)}+(1-n)(2g^{(2)})-2g^{kl}g^{(2)}_{kl}g-2Ric(g)=0.
\ee
Taking trace with respect to $g$, which we will denote as $\Tr_{g}$, we
get
$$
R_{g}=(2-2n)\Tr_{g} g^{(2)},
$$or,
\be\label{eqn6}
J_{g}={R_{g}\over (2n-2)}=-\Tr_{g}g^{(2)}. \ee
Combine (\ref{eqn5}) and (\ref{eqn6}), we get
\be\label{eqn66}
g^{(2)}={1\over (2-n)}(Ric_{g}-J_{g})=-P_{g}. \ee
and
$$ v^{(2)}={{-1}\over 2} J_{g}. $$

We now apply the same method to compute $v^{(4)}$. Studying the
coefficient of $r^{3}$ in (\ref{eqn4}), we have \be\label{eqn7} 12
g^{(4)}+(1-n) (4g^{(4)}) - \alpha -\beta+\gamma-2\delta=0, \ee
where $\alpha,\beta,\gamma,\delta$ are the coefficients of $r^{3}$
in $\Tr [h^{-1}h'] h,\ r h^{'}h^{-1}h',\ {r\over 2}h^{-1}h'h'$ and
$r\  Ric(h)$, respectively.

We now compute these coefficients. First notice that
\ba h^{-1}&=&[(\rm{Id}+B) g]^{-1}=g^{-1}
(\rm{Id}-B+B^{2}-B^{3}+\cdots),\nonumber \\
&=&g^{-1}(1-r^{2}C_{2}-r^{4}C_{4}+r^{4}C_{2}^{2}+o(r^{4})).\label{eq7}
\ea
Combing with (\ref{eqn4}), we get
\be\label{alpha}
\alpha=-2(\Tr\ C_{2}^{2})g+\Tr_{g}
[(4g^{(4)})]g+\Tr_{g}[2g^{(2)}]g^{(2)}.\ee

Similarly, we can get,
\be\label{beta} \beta=2g^{(2)}\ g^{-1}\ (2g^{(2)}),\ee
\be\label{gamma}\gamma={1\over 2}\Tr_{g}[2g^{(2)}](2g^{(2)}),\ee
Thus, by (\ref{eqn7}), we get
\ba
0= 12 g^{(4)}+(1-n) (4g^{(4)})+2(\Tr\ C_{2}^{2})g-\Tr_{g} [(4g^{(4)})]g
\nonumber \\ -\Tr_{g}[2g^{(2)}]g^{(2)}-4g^{(2)}\ g^{-1}\
g^{(2)}+2\Tr_{g}[g^{(2)}](g^{(2)})-2\delta
\ea

Taking trace with respect to $g$, we get
$$
(16-8n)\Tr_{g}g^{(4)}+(2n-4)(\Tr \ C_{2}^{2})-2\Tr_{g}\delta=0;$$ or,
equivalently we get the fomula
\be\label{eqn71}
[ (8-4n)\Tr\ C_{4}+(n-2)(\Tr \ C_{2}^{2})]-\Tr_{g}\delta=0.
\ee
Regarding to the term involving $\delta$, we will prove the following

\begin{lemma}\label{Ricci}
$$\Tr_{g}\delta=(4-4n)\Tr\  C_{4}+(n-1) \Tr \ C_{2}^{2}.$$
\end{lemma}
         
This lemma can be verified by relating the scalar curvature of $\bar g $
to that of the scalar curvature of $g_{+}$ and $h (r, .)$. The computation
is tedious but relatively routine. We will skip the detail here.

Combine the formula in \ref{eqn71} and~\ref{Ricci}, we have
\be\label{C4}\Tr\ C_{4}=\Tr \ C_{2}^{2}.\ee

By Lemma~\ref{v2v4} and (\ref{C4}), we get
$$v^{(4)}={1\over 2}[\Tr\  C_{4}+\sigma_{2}(C_{2})-({1\over 4})(\Tr
C^{2}_{2})]={1\over 2}(\sigma_{2}(C_{2})-{1\over 2}\sigma_{2}(C_{2}))=
{1\over 4}\sigma_{2}(C_{2}).$$ Noticing that $C_{2}=-P_{g}$, where
$P_{g}$ is the Schouten tensor, we have established the following
\begin{theorem}
For $n>4$, we have
$$v^{(4)}(g)={1\over 4}\sigma_{2}(g).$$
\end{theorem}

\vskip .2in

\section{Variational property of ${ v}^{(6)}$}

In this section, we study the properties of
$$
{ v}^{(6)}(g) \, = \, - \frac {1}{8} \left[\sigma_{3}(g) \, + \,  {1\over{3(n-4)}}(P_g)^{ij}{(B_g)}_{ij}\right].
$$ 

We give a direct proof of the following special case of Theorem~\ref{main}.
\begin{theorem}\label{v6'}

 For any metric $h$ on $M$ and $6\leq n=\dim M$, define a functional
 \begin{equation}
 {\cal F}_{3}(h)=\int_{M}{ v}^{(6)}(h)dvol_{h}/
(\int_{M}dvol_{h})^{{n-6}\over n},\nonumber
 \end{equation}
 then ${\cal F}_{3}$ is variational within the conformal class; i.e., the critical metric in $[h]$ satisfies the equation
 \begin{equation}
 { v}^{(6)}=constant.
 \end{equation}
 For $n=6$, we have that $F_{3}(h)$ is constant in the conformal class $[h]$.

\end{theorem}

To prove the theorem, we first recall some basic conformal transformation law for the curvature
invariants involved.

\begin{lemma}
For a fixed smooth function $\phi$ defined on $M$, we consider two conformal equivalent metrics, $g$ and  $g_{\phi}= e^{2\phi} g$  on $M$. Then we have, under a local coordinate system,

\begin{equation}
P{(g_{\phi})}_{ij}=P_{ij}-\phi_{ij}+\phi_{i}\phi_{j}-\frac{1}{2}|\nabla
\phi|^{2}g_{ij},\label{P}
\end{equation}
\begin{equation}
B(g_{\phi})_{ij}=e^{-2\phi}[B_{ij}-(n-4)(C_{ikj}+C_{jki}){\phi}^{k}-(n-4)W_{kijl}{\phi}^{k}{\phi}^{l}],\label{B}
\end{equation}
where $W$ and $C$ are the Weyl tensor and Cotten tensor of $g$, respectively.
\end{lemma}

We now consider a family of conformal metrics on $M$, $g_{t}=e^{2t\phi}g$ and denote $D= {d\over{dt}}|_{t=0}$ as before.

Now we can compute $D [{\cal F}_{3}(g_{t})]$. We separate the computation in two steps.

First,
\ba
& &D \ \int {B_{ij}(g_{t})P^{ij}(g_{t})dvol (g_{t} )}\nonumber \\
&=&\int {\{D[B_{ij}(g_{t})]P^{ij}+B_{ij} D[P^{ij}(g)t] \}dvol +B_{ij}P^{ij}D[dvol(g_{t})]\} }\nonumber \\
&=&\int{  [-2 \phi B_{ij}P^{ij} -2(n-4)C_{ijk}P^{ij} {\phi}^{k}-B_{ij} {\phi}^{ij}+(n-4) \phi B_{ij}P^{ij}]dvol }\nonumber \\
&=&\int{  [-2 \phi B_{ij}P^{ij} -2(n-4)C_{ijk}P^{ij} {\phi}^{k}+\nabla^{i}B_{ij} {\phi}^{j}+(n-4) \phi B_{ij}P^{ij}]dvol }\nonumber \\ 
&=&\int{  [ -2(n-4)C_{ijk}P^{ij} {\phi}^{k}-(n-4)C_{ijk}P^{ij} {\phi}^{k}+(n-6) \phi B_{ij}P^{ij}]dvol }\nonumber \\
&=&\int{  [  -3(n-4)C_{ijk}P^{ij} {\phi}^{k} +(n-6) \phi B_{ij}P^{ij}]dvol },\nonumber
\ea

Second, define the Newton tensor as
$$T^{ij}=\sigma_{2}(g) g^{ij} -\sigma_{1}(g) P^{ij}+ (P^{2})^{ij}.$$
Then we have

\ba
&&D \ [\int {\sigma_{3}(g_{t})dvol(g_{t})}]\nonumber \\
&=&\int {[\sigma_{3}(g)(n-6) \phi -T^{ij} {\phi}_{ij}] dvol}\nonumber \\
&=&\int {[\sigma_{3}(g)(n-6) \phi +T^{ij}_{,j}{\phi}_{i}] dvol}.
\ea

By the Bianchi identity, we have
\ba
P^{ij}_{,j}&=&\nabla^{i}J;\nonumber \\
(P^{2})^{ij}_{,j}&=&(P^{k}_{i}P^{kj})_{,j}\nonumber \\
&=&P^{i}_{k}J^{,k}+P^{i}_{k,j}P^{kj}\nonumber \\
&=&P^{i}_{k}J^{,k}+ g^{is}C_{ksj}P^{kj}+P_{kj}P_{kj}^{,i}\nonumber \\
&=&P^{i}_{k}J^{,k}+ C^{kij}P_{kj}+ \frac {1}{2} \nabla^{i}[\Tr_{g}(P^{2})]\nonumber \\
&=&P^{ik}J_{,k}+ C^{kij}P_{kj}+ \frac {1}{2} \nabla^{i}[\Tr_{g}(P^{2})].
\ea
Thus, we have
\ba
& &\int {T^{ij}_{,j}{\phi}_{i} dvol}\nonumber \\
&=&\int { \{\sigma_{2}^{,i} -J_{,j}P^{ij}-JJ^{,i}+P^{ik}J_{,k} +C^{kij}P_{kj}
+ \frac{1}{2} \nabla^{i}[\Tr_{g}(P^{2})]\} \phi_{i}dvol}\nonumber \\
&=&\int{ [C^{kij}P_{kj}{\phi}_{i}+\nabla^{i}[\sigma_{2}+{1\over 2}(\Tr_{g}P^{2}-J^{2})]{\phi}_{i} dvol}\nonumber \\
&=&\int{ C^{kij}P_{kj}{\phi}_{i}dvol}.
\ea

Finally, we can combine these to get
\ba
&&D{\cal F}_{3}(g_{t})\nonumber \\
&=&
- \frac 18 [\, { D \ (\int
\sigma_{3}(g_{t})dvol(g_{t}))+{1\over{3(n-4)}}D \ \int
{B_{ij}(g_{t})P^{ij}(g_{t})dvol (g_{t} ) }}]\nonumber \\
&=& (n-6)\int { v}^{(6)} \phi dvol.
\ea
Theorem~\ref{v6'} then follows easily.

\begin{remark}
From (\ref{P}) and (\ref{B}), we see that the equation
\begin{equation} \label{var-v}
{ v}^{(6)} (g_w) ={\rm constant} .\end{equation}
is a second
order PDE in terms of the conformal factor $w$. This is in
analogue of $\sigma_{k}(g_w) ={\rm constant} $ equation which has been
intensively studied in recent years. It
remains to see under what conditions can the PDE ({\ref{var-v}})
be solved for metrics in a fixed conformal class and if
the sign of the integral $\int { v}^{(6)} (g) dv_g $ plays
some role and carries
geometric information as in the case 
for $\int {v}^{(4)}(g)
dv_g $ on manifolds of dimension 3 and 4 (cf. [GV01, CGY02]). The authors wish to report some further study of this problem
in the future.

\end{remark}

\end{document}